\documentclass[final,1p,times]{elsarticle}

\usepackage{graphicx}
\usepackage{color}
\newtheorem{thm}{Theorem}

\usepackage{amssymb,amsmath}

\usepackage{multirow}
\usepackage{stackrel}

\newdefinition{rmk}{Remark}

\DeclareMathOperator{\dd}{d}

\journal{XXX}
\usepackage[all]{xy}

\begin{document}

\begin{frontmatter}



\title{Solving the Random Pielou Logistic Equation with the Random Variable  Transformation Technique: Theory and Applications}


\author[IMM]{J.-C. Cort\'es\corref{cor1}}
\ead{jccortes@imm.upv.es}

\author[deusto1,deusto2]{A. Navarro-Quiles}
\ead{annaqui@doctor.upv.es}

\author[IMM]{J.-V. Romero}
\ead{jvromero@imm.upv.es}

\author[IMM]{M.-D. Rosell\'{o}}
\ead{drosello@imm.upv.es}

\address[IMM]{Instituto Universitario de Matem\'{a}tica Multidisciplinar,\\
Universitat Polit\`{e}cnica de Val\`{e}ncia,\\
Camino de Vera s/n, 46022, Valencia, Spain}

\address[deusto1]{DeustoTech, University of Deusto, 48007 Bilbao,\\
 Basque Country, Spain}
\address[deusto2]{ Facultad de Ingenieria, Universidad de Deusto, \\
Avda.Universidades, 24, 48007, Bilbao, Basque Country, Spain.}

\cortext[cor1]{Corresponding author. Phone number: +34--963877000 (ext.~88289)}

\begin{abstract}
The study of the dynamics of the size of a population via mathematical modelling is a problem of interest and widely studied. Traditionally, continuous deterministic methods based on differential equations have been used to deal with this problem. However discrete versions of some models are also available and sometimes more adequate. In this paper, we randomize  the Pielou logistic equation in order to include the inherent uncertainty in  modelling. Taking advantage of the method of transformation of random variables, we provide a full probabilistic description to  the randomized Pielou logistic model via the computation of the probability density functions of the solution stochastic process, the steady state and the time until a certain level of population is reached. The theoretical results are illustrated by means of two examples, the first one consists of a numerical experiment and the second one shows an application to study the diffusion of a technology using real data.
\end{abstract}


\end{frontmatter}

\section{Introduction}
\label{sec1}

Some interesting problems in population dynamics are to model the changes 
in the size and to quantify the composition of the population along the time\cite{molles2007ecology}.
Mathematical models of populations can be used to accurately describe changes in a population and, more importantly, to predict future changes.

Population growth can be described by several mathematical models. The  simplest model is the Malthusian one.  It can predict an exponential increase in
the population with time, which is unrealistic for long time. But it can be an adequate model in short-time studies. A typical application of Malthusian model is the study of  the evolution of bacteria in a laboratory.

A more realistic and popular continous model of  population growth is the logistic o Verhulst-Pearl equation
\begin{equation}\label{logistic}
x'(t)=x(t)(\alpha-\beta x(t)), \quad \alpha,\beta>0,
\end{equation}
where $x(t)$ is the size of the population at time $t$, $\alpha$ is the rate of population growth if the resources were unlimited and the individuals did not affect  one another, and $-\beta x^2(t)$ represents the negative effect on the growth due to crowdedness and limited resources. This model was introduced by Pierre Verhulst \cite{ver,ver2}.  Although originally this equation was introduced to model populations, it has also been used to model different problems \cite{KWASNICKI201350}.

Most models that describe population dynamics are continuous. But many times only data is available for discrete times, so it is interesting to have the solution for those moments. In this case, a difference equation instead of a differential equation can be more appropriate to model this situation.  
For the sake of clarity, as time variable can be considered discrete or continuous, depending upon the context, hereinafter we will called period in the former case and time in the latter.

The discrete version of the logistic equation is known as Pielou logistic equation and it
is stated as
\begin{equation} \label{probdet}
x_{n+1}=\frac{a x_n}{1 + b x_n},
\end{equation}
where $a>1$ and $b>0 $\cite{elaydi2005,pielou1969,pielou1974}.
By letting  $z_n=1/x_n$, the nonlinear difference Equation \eqref{probdet} is transformed into the following linear difference equation
\begin{equation} \label{lineardet}
 z_{n+1} = \frac{1}{a} z_n + \frac{b}{a}.
\end{equation}
Taking the initial condition $z_0=1/c$, the solution of Equation \eqref{lineardet} is
\[ 
 z_n =
 \left\{
   \begin{array}{llll}
   \displaystyle 
     \left[ \frac{1}{c} - \frac{b}{a-1} \right] a^{-n} + \frac{b}{a-1}, &&& \text{if } a \neq 1, \\ \\
   \displaystyle 
      \frac{1}{c} + b n, &&& \text{if } a = 1.
   \end{array}
 \right.
\]

Thus, the solution of Pielou logistic equation  \eqref{probdet} is given by
\begin{equation}  \label{soldet}
 x_n =
 \left\{
   \begin{array}{llll}
   \displaystyle 
     \frac{a^n (a -1)}{b a^n + \frac{1}{c} (a -1)-b}, &&& \text{if } a \neq 1, \\
     \\
   \displaystyle 
      \frac{1}{\frac{1}{c}+b n} &&& \text{if } a = 1.
   \end{array}
 \right.
\end{equation}

In real problems coefficients and initial conditions are not usually known exactly. 
This may be due to measurement errors or the
inherent complexity associated to their own nature. So, it  seems more realistic to consider that parameters and initial conditions are RVs instead of deterministic values. Notice that hereinafter capital letters are used to denote a RV.

There are some recent interesting contributions concerning continuous random modelling in population dynamics\cite{Chen2011,dorini2016}. 
Discrete random models have not been widely studied in the literature. Recently, authors of this paper have made several contributions related to Markov models\cite{Cortes2017estroke,Cortes2018IJCM}  and linear difference equations\cite{CORTES2017150}.

The main objective of this paper is to construct a randomized version of the Pielou logistic equation. To the best of our knowledge, this problem
has not been considered in the extant literature yet, but randomizing parameters of a model is a
technique used in other contexts.

Solving a random difference  equation means not only to calculate the exact solution of the stochastic process but also its main statistical functions. 
The computation of the 1-PDF, $f_1^X(x;n)$, allows us to have a complete statistical description of the solution. From the 1-PDF we can calculate easily the mean, variance and other higher statistical moments, by means of the following expressions, respectively
\begin{equation}
 \mu_X(n)=\mathbb{E}[X_n]=\int_{-\infty}^{\infty} x\,f_1^X(x;n)\mathtt{d}x,
 \label{expectation}
\end{equation}

\[
 \sigma^2_X(n)=\mathbb{V}[X_n]=\int_{-\infty}^{\infty} x^2\,f_1^X(x;n)\mathtt{d}x -\left( \mu_X(n) \right)^2,
\]

\[
\mathbb{E}[(X_n)^k]=\int_{-\infty}^{\infty} x^k\,f_1^X(x;n)\mathtt{d}x,\quad k=0,1,2,\ldots
\]

Furthermore significant information such as the probability of the solution lies within a set of interest can be determined from the 1-PDF
\[
\mathbb{P}[a \leq X_n \leq b]=\int_{a}^{b} f_1^X(x;n)\mathtt{d}x.
\]

This improves the computation of rough bounds, like the one derived via Chebyshev's inequality\cite{Soong}
\[
\mathbb{P}\left[ | X_n -\mu_{X}(n) | \geq  \lambda \right]
\leq 
\frac{(\sigma_{X}(n))^2}{\lambda^2} \,,\quad \lambda>0,
\]
usually applied in practice.

RVT method is a powerful technique  that has been recently used by the authors to compute the 1-PDF of the solution of some differential and difference equations \cite{CORTES2018190,CASABAN2017396,CORTES2017225,CORTES2017150}.
The RVT technique permits to compute the PDF of a RV which results from mapping another RV whose PDF is  known. The multidimensional version of the RVT technique is stated in Theorem \ref{R.V.T.method2}.

\begin{thm}[RVT Multidimensional version\cite{Soong,CASABAN2017396} \label{R.V.T.method2}]
Let $\mathbf{U}=(U_1,\ldots,U_m)^{\top}$ and $\mathbf{V}=(V_1,\ldots,V_m)^{\top}$ be two $m$-dimensional absolutely continuous random vectors. Let $\mathbf{r}: \mathbb{R}^m \rightarrow \mathbb{R}^m$ be a one-to-one deterministic transformation of $\mathbf{U}$ into $\mathbf{V}$, i.e., $\mathbf{V}=\mathbf{r}(\mathbf{U})$. Assume that $\mathbf{r}$ is continuous in $\mathbf{U}$ and has continuous partial derivatives with respect to $\mathbf{U}$. Then, if $f_{\mathbf{U}}(\mathbf{u})$ denotes the joint probability density function of vector $\mathbf{U}$, and $\mathbf{s}=\mathbf{r}^{-1}=(s_1(v_1,\ldots,v_m),\ldots,s_m(v_1,\ldots,v_m))^{\top}$ represents the inverse mapping of $\mathbf{r}=(r_1(u_1,\ldots,u_m),\ldots,r_m(u_1,\ldots,u_m))^{\top}$, the joint probability density function of vector $\mathbf{V}$ is given by
\begin{equation}\label{RVTformula}
f_{\mathbf{V}}(\mathbf{v})=f_{\mathbf{U}}\left(\mathbf{s}(\mathbf{v})\right) \left| J_m \right|,
\end{equation}
where $\left| J_m \right|$ is the absolute value of the Jacobian, which is defined by
\begin{equation}\label{jacobiano_general}
J_m=\det \left( \frac{\partial \mathbf{s}^\top}{\partial \mathbf{v}}\right)
=
\det
\left(
\begin{array}{ccc}
\dfrac{\partial s_1(v_1,\ldots, v_m)}{\partial v_1} & \cdots & \dfrac{\partial s_m(v_1,\ldots, v_m)}{\partial v_1}\\
\vdots & \ddots & \vdots\\
\dfrac{\partial s_1(v_1,\ldots, v_m)}{\partial v_m} & \cdots & \dfrac{\partial s_m(v_1,\ldots, v_m)}{\partial v_m}\\
\end{array}
\right)
.
\end{equation}
\end{thm}

\begin{rmk}
In the context of solving random difference equations, when Theoerem \ref{R.V.T.method2} is applied the choice of random vector $\mathbf{V}$ (that defines the mapping $\mathbf{r}$) often can be made in different ways. The convenience of choosing one mapping against another depends heavily  upon the feasibility   of obtaining the inverse mapping $\mathbf{s}$ as well as the involved computations to obtain the 1-PDF of the solution stochastic process. Anyway, the 1-PDF obtained is equivalent independently of the mapping initialy chosen. This issue will be illustrated later (see Remark \ref{rem2}).
\end{rmk}

This paper is organized as follows. In Section \ref{sec2},  the Pielou logistic equation \eqref{probdet} is randomized.  Then, the solution, the steady state and the time until a given population level is achieved are determined. Also a full probabilistic description of these magnitudes is provided via their 1-PDF and PDFs, respectively. In Section \ref{sec3}, the theoretical findings are illustrated by examples.
Finally, conclusions are drawn in Section \ref{sec4}.

\section{Randomized Pielou logistic equation}
\label{sec2}

Random Pielou logistic equation can be written as
\begin{equation} \label{Pielourand}
\left\{
\begin{array}{l}
   \displaystyle X_{n+1}= \frac{A X_n}{1+B X_n}, \qquad n=0,1,2,\ldots  \\ 
   X_{0}=C,
\end{array}
\right.
\end{equation}
where all the input parameters $A$,$B$ and $C$ are assumed to be absolutely continuous RVs defined on a common complete probability space 
($\Omega, \mathcal{F},\mathbb{P}$). As a natural extension of its deterministic counterpart\citep{elaydi2005}, we assume 
$\mathbb{P}\left[ \{ \omega \in \Omega:\, A(\omega)>1 \} \right]=1$,
$\mathbb{P}\left[ \{ \omega \in \Omega:\, B(\omega)>0 \} \right]=1$
and $\mathbb{P}\left[ \{ \omega \in \Omega:\, C(\omega)>0 \} \right]=1$. For the sake of generality, hereinafter we will assume that $A,B,C$ are dependent RVs whose joint PDF is $f_{C,A,B}(c,a,b)$.

The main goal of this section is to obtain the 1-PDF of the solution of the random Pielou logistic equation \eqref{Pielourand}, say $f_1^{X}(x;n)$.
In this section we will also determine the PDF of another interesting quantities in dealing with the Pielou equation. Specifically, the steady state and the time until a given proportion is reached will be studied from a probabilistic standpoint.

\subsection{1-PDF of the solution of the randomized Pielou equation}
\label{subsec21}

Inspired in the deterministic theory \eqref{probdet}--\eqref{soldet}, by introducing the change of variable $Z_n=1/X_n$ in Equation \eqref{Pielourand}, it is linearized and then solved. 
As $A$ is an absolutely continuous RV, then  $\mathbb{P}\left[ \{ \omega \in \Omega:\, A(\omega)=1 \} \right]=0$, for all event $\omega\in \Omega$. As a consequence,  we obtain that  the solution of random Pielou logistic problem \eqref{Pielourand} is given by
\begin{equation} \label{sol}
X_n
= \frac{ A^n (A-1)}{ B A^n+ \frac{1}{C} (A-1)-  B}
,\qquad n=0,1,\ldots
\end{equation}
and this solution is well-defined from a probabilistic point of view.

In order to compute the 1-PDF of \eqref{sol}, we can define the transformation  $\mathbf{r}:\mathbb{R}^3\longrightarrow \mathbb{R}^3$, 
\[
\begin{array}{ccccl}
y_1 & = & r_1\left( c,a,b \right) &=& \displaystyle \frac{a^n (a-1)}{ba^n+ \frac{1}{c} (a-1)-b},\\
y_2 & = & r_2\left( c,a,b \right) &=& a,\\
y_3 & = & r_3\left( c,a,b \right) &=& b,
\end{array}
\]
whose inverse mapping $\mathbf{s}=\mathbf{r}^{-1}$ is  given by
\[
\begin{array}{lllll}
c & = & s_1 \left(y_1,y_2,y_3  \right) &=& \displaystyle
\frac{y_1 \left(y_2 -1 \right)}{y_2^n \left( y_2 -1 \right) - y_3 y_1 \left( y_2^n -1 \right)},\\
a & = & s_2 \left(y_1,y_2,y_3 \right) & =&y_2 ,\\
b & = & s_3 \left(y_1,y_2,y_3 \right) &=& y_3, 
\end{array}
\]
and the absolute value of the Jacobian of the inverse mapping, $\mathbf{s}$, is
\[
\left | J_3 \right | = \left |  \frac{\partial s_1}{\partial y_1}\right|=
\left |  \frac{\left( y_2 -1 \right)^2 y_2^n}{\left( y_2^n \left( y_2 -1 \right) - y_3 y_1 \left( y_2^n -1 \right)  \right)^2}\right|,
\]
which is different from zero if $y_2=a>1$, which holds by hypothesis.

Applying  RVT technique (Theorem \ref{R.V.T.method2}) for an arbitrary but fixed value of $n$, the PDF of the random vector $(Y_1,Y_2,Y_3)$ defined by mapping $\mathbf{r}$ is given by
\begin{equation} \label{PDFconj}
f_{Y_1,Y_2,Y_3}\left( y_1,y_2,y_3\right) =
f_{C,A,B}
\left(
\frac{y_1 \left(y_2 - 1 \right)}{y_2^n \left( y_2 -1 \right) - y_3 y_1 \left( y_2^n -1 \right)},
y_2,y_3
\right)
\left |  \frac{\left( y_2 -1 \right)^2 y_2^n}{\left( y_2^n \left( y_2 -1 \right) - y_3 y_1 \left( y_2^n -1 \right)  \right)^2}\right|.
\end{equation}
Finally, the 1-PDF of the solution of the randomized Pielou logistic model, $X_n$ is obtained by marginalizing expression \eqref{PDFconj} with respect to $A$ and $B$, being $n$ arbitrary,
\begin{equation}
\label{f1xn}
f_1^{X}(x;n) = 
\iint_{\mathcal{D}(A,B)} 
f_{C,A,B}
\left(
\frac{x \left(a -1 \right)}{a^n \left( a -1 \right) - b x \left( a^n -1 \right)},
a,b
\right)
\left|  \frac{\left( a -1 \right)^2 a^n}{\left( a^n \left( a -1 \right) - b x \left( a^n -1 \right)  \right)^2} \right| \dd\!a \dd\!b,
\end{equation}
where $\mathcal{D}(A,B)$  stands for the domain of the random vector $(A,B)$.

\begin{rmk}\label{rem2}
Notice that another possibility to choose an adequate mapping  $\mathbf{r}$ when applying Theorem \ref{R.V.T.method2} would be
\[
\begin{array}{ccccl}
y_1 & = & r_1\left( c,a,b \right) &=& \displaystyle \frac{a^n (a-1)}{ba^n+ \frac{1}{c} (a-1)-b},\\
y_2 & = & r_2\left( c,a,b \right) &=& a,\\
y_3 & = & r_3\left( c,a,b \right) &=& c.
\end{array}
\]
This choice would lead to an equivalent expression to  $f_1^{X}(x;n)$ given in \eqref{f1xn}.
\end{rmk}

\subsection{PDF of the steady state of the randomized Pielou equation}
\label{subsec22}

A key magnitude in population dynamics is to determine its long term behaviour. Mathematically it can be obtained calculating the limit as $n$ tends to infinite of the solution $X_n$. So, the steady state of the randomized Pielou equation, $X_{\infty}$, is
\begin{equation}\label{estacionariorand}
X_\infty= \lim_{n\to \infty} X_n=
\lim_{n\to \infty}
\frac{ {A}^n ({A}-1)}{ {B} {A}^n+ \frac{1}{{C}} ({A}-1)-   {B}}
=
\frac{{A}-1}{{B}}.
\end{equation}
Notice that,  as  $\mathbb{P}\left[ \{ \omega \in \Omega:\, A(\omega)>1 \} \right]=1$  and $\mathbb{P}\left[ \{ \omega \in \Omega:\, B(\omega)>0 \} \right]=1$    for all event $\omega\in \Omega$ , $X_\infty$ is a well-defined positive RV as is required for the size of a population.

By applying again the RVT technique (Theorem \ref{R.V.T.method2}), considering the mapping  $\mathbf{r}:\mathbb{R}^3\longrightarrow \mathbb{R}^3$, 
\[
\begin{array}{ccccl}
y_1 & = & r_1\left( c,a,b \right) &=& \displaystyle \frac{a-1}{b},\\
y_2 & = & r_2\left( c,a,b \right) &=& a,\\
y_3 & = & r_3\left( c,a,b \right) &=& b,
\end{array}
\]
and, after some technical computations, we obtain the PDF of the RV steady state, 
\begin{equation}
f_{X_{\infty}}(x)=
\iint_{\mathcal{D}(C,B)}
f_{C,A,B}
\left( c, x b+1,b  \right)
\left| b \right| \dd\!c \dd\!b,
\label{pdfinf}
\end{equation}
where $\mathcal{D}(C,B)$  stands for the domain of the random vector $(C,B)$.

\subsection{PDF of time until a given population size is reached}
\label{subsec23}

Another interesting question in dealing with discrete population models is
to determine the distribution of the time  where the size of the population reaches a certain specific value, $\hat{X}$,
\begin{equation}\label{xgorro}
\hat{X}=X_n
= \frac{ {A}^n ({A}-1)}{ {B} {A}^n+ \frac{1}{{C}} ({A}-1)-   {B}}.
\end{equation}
Again the RVT technique is very useful to answer this interesting question. In order to determine the PDF of the time where the size of the population reaches a certain specific value, say $f_N(n)$, we first isolate $n$ from expression \eqref{xgorro}, and as the obtained expression is a RV, we use capital letter notation to denote it, $N$, 
\[
N=\dfrac{\log \left( \frac{\hat{X} ({C} {B} - {A} +1)}{{C} ({B} \hat{X} - {A} +1)}  \right)}{\log({A})}.
\]

Defining  an appropriate mapping, for example, $\mathbf{r}:\mathbb{R}^3\longrightarrow \mathbb{R}^3$, 
\[
\begin{array}{ccccl}
y_1 & = & r_1\left( c,a,b \right) &=& \displaystyle \frac{\log \left( \frac{\hat{x} ({c} {b} - {a} +1)}{{c} ({b} \hat{x} - {a} +1)}  \right)}{\log({a})},\\
y_2 & = & r_2\left( c,a,b \right) &=& a,\\
y_3 & = & r_3\left( c,a,b \right) &=& b,
\end{array}
\]
and applying the RVT method (Theorem \ref{R.V.T.method2}), 
 one can obtain the PDF of $N$,
\begin{equation}
f_{N }(n) = \iint_{\mathcal{D}(A,B)}
f_{C,A,B}
\left( \dfrac{\hat{x}(1-a)}{\hat{x}b(a^n-1)+a^n(1-a)}, a,b  \right)
 \dfrac{\left|\hat{x}(a-1)a^n(1-a+\hat{x} b)\log(a)\right|}{\left( \hat{x} b (a^n-1)+a^n(1-a) \right)^2}  \dd\!a \dd\!b,
\label{pdftiempo}
\end{equation}
where $\mathcal{D}(A,B)$  stands for the domain of the random vector $(A,B)$.

\section{Examples}
\label{sec3}

In this section we will present two examples. 
In Subsection \ref{subsec31}, we show a numerical example aimed to illustrate the theoretical results obtained in Section \ref{sec2}. In Subsection \ref{subsec32}, a second example is presented concerning modelling. Although Pielou equation is usually applied to model population dynamics, it can model other interesting situations as the diffusion of a technology. In particular, we describe the dynamics of the number of mobile lines in Spain during the range of years 1999--2015 using real data by means of the randomized Pielou equation. We assume that $A$, $B$ and $C$ are independent RVs.

\subsection{Numerical example}
\label{subsec31}

We consider the randomized Pielou model \eqref{Pielourand} and we assume the following probability independent distributions for the inputs. For the initial condition we choose a truncated normal distribution, ${C}\sim \text{N}_{]0,1[}(0.5;0.05)$, and for $A$ and $B$ uniform and beta distributions ${A}\sim \text{Un}([1.1,2])$ and ${B}\sim \text{Be}(2;3)$, respectively. 
We have chosen these distributions to illustrate the capability of our approach to deal with different probability distributions, although other distribution can be used too.
All the computations have been carried out using the software Mathematica\textsuperscript{\textregistered}\cite{Mathematica}.

\begin{figure}[htp]
  \centerline{\includegraphics[width=0.6\textwidth]{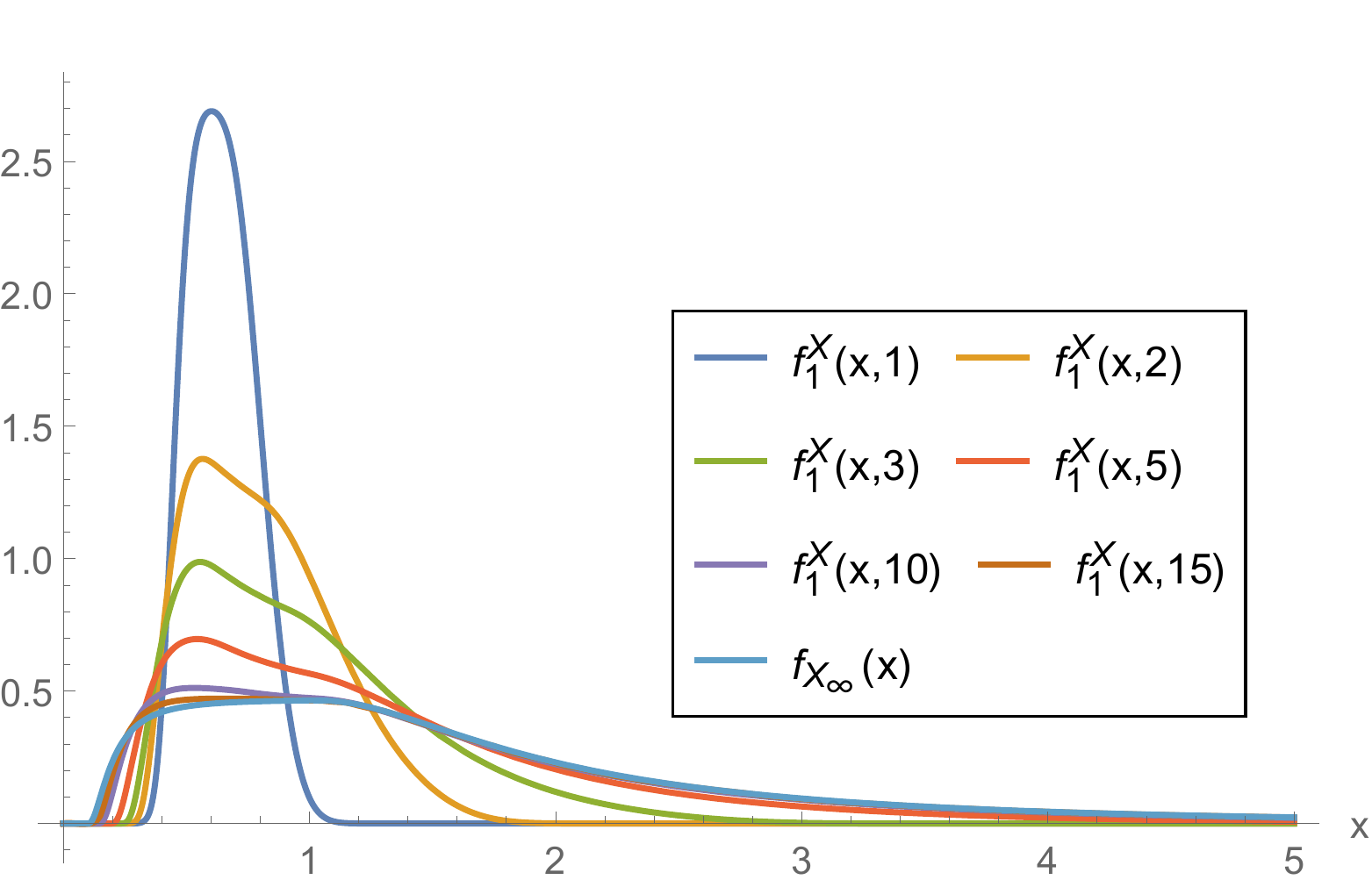}}
\caption{1-PDF of the solution stochastic  process $X_n$ with $n\in\{1,2,3,5,10,15\}$ and PDF of the equilibrium RV $X_{\infty}$. Example \ref{subsec31}.}
\label{fig1}
\end{figure}

\begin{figure}[htp]
\centerline{
 \includegraphics[width=0.6\textwidth]{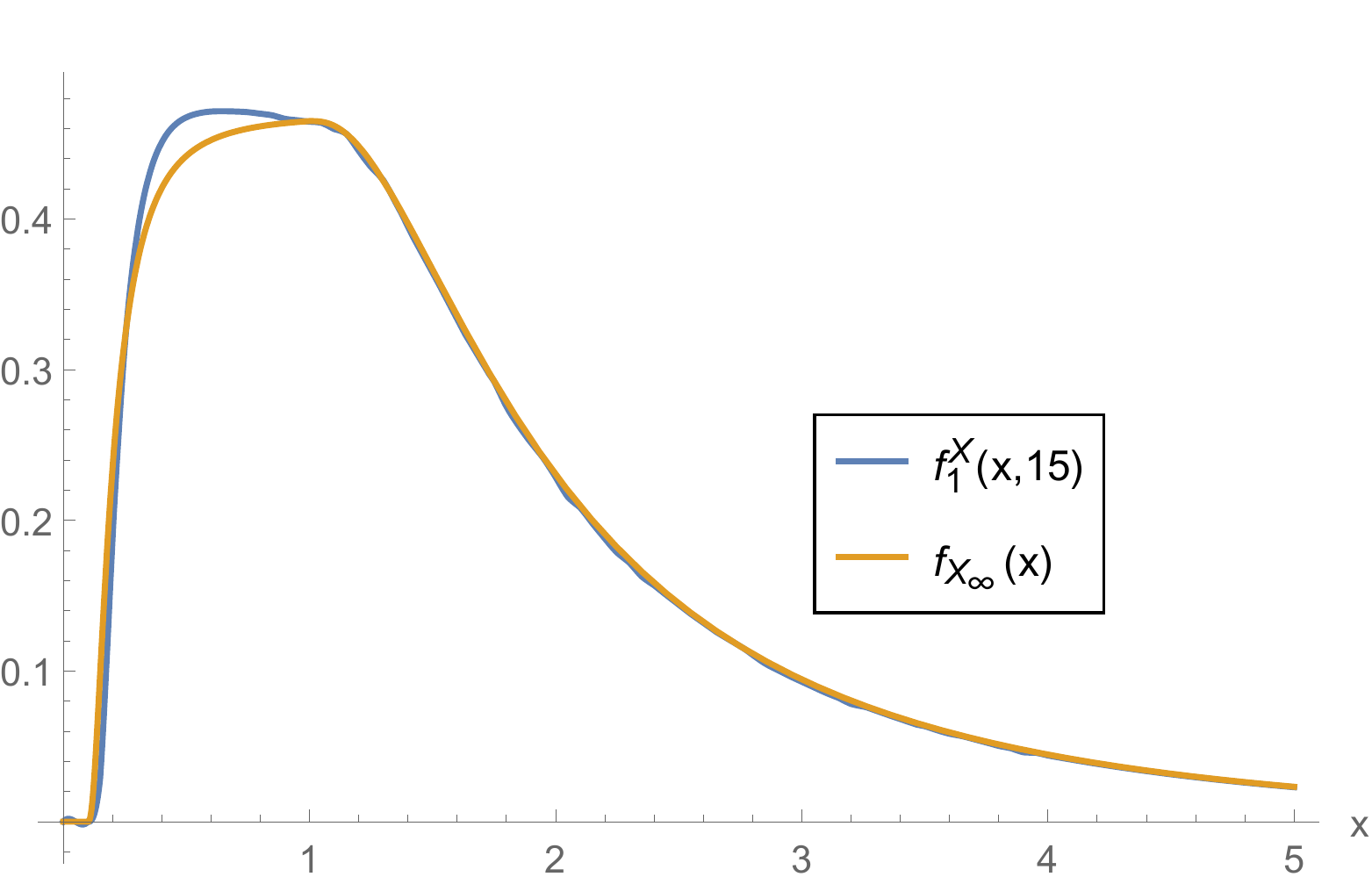}
}
\caption{1-PDF of the solution stochastic  process in $n=15$, $X_{15}$, and PDF of the steady state $X_{\infty}$. Example \ref{subsec31}.}
\label{fig2}
\end{figure}

We have calculated the 1-PDF of the solution stochastic  process $X_n$, $f_1^{X}(x;n)$,  given by \eqref{f1xn}. Also we have determined the PDF of the steady state,  $f_{X_\infty}(x)$, given by \eqref{pdfinf}.  In Figure \ref{fig1} it is plotted $f_1^{X}(x;n)$ for different values of $n$ and 
$f_{X_\infty}(x)$. We can observe that $f_1^{X}(x;n)$ tends to $f_{X_\infty}(x)$ as $n$ increases. For the sake of clarity, in Figure \ref{fig2} it is plotted $f_1^{X}(x;15)$  and $f_{X_\infty}(x)$, and we can observe that both PDFs are similar.

\begin{figure}[htp]
\centerline{
\includegraphics[width=0.45\textwidth]{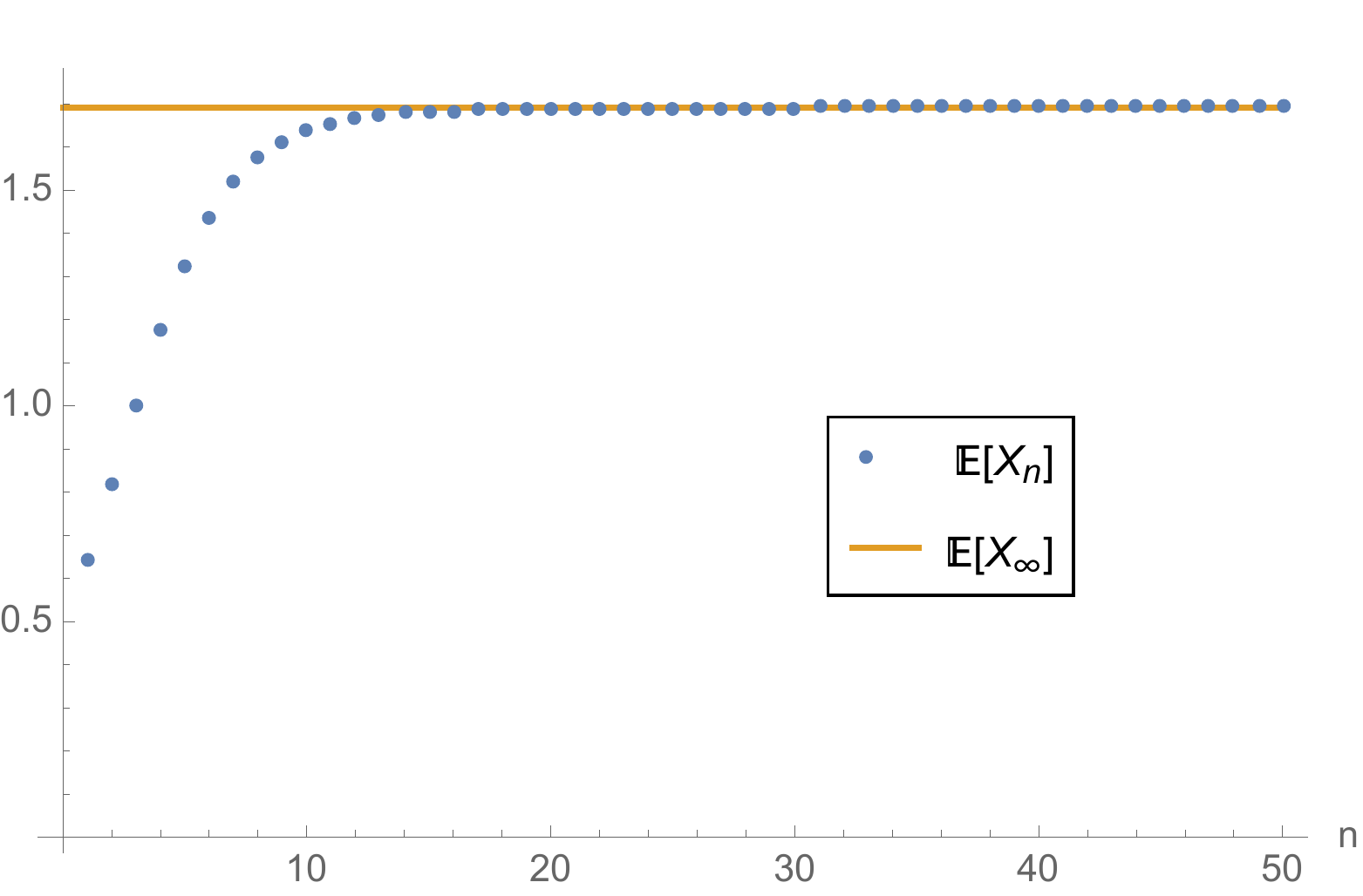}
\includegraphics[width=0.45\textwidth]{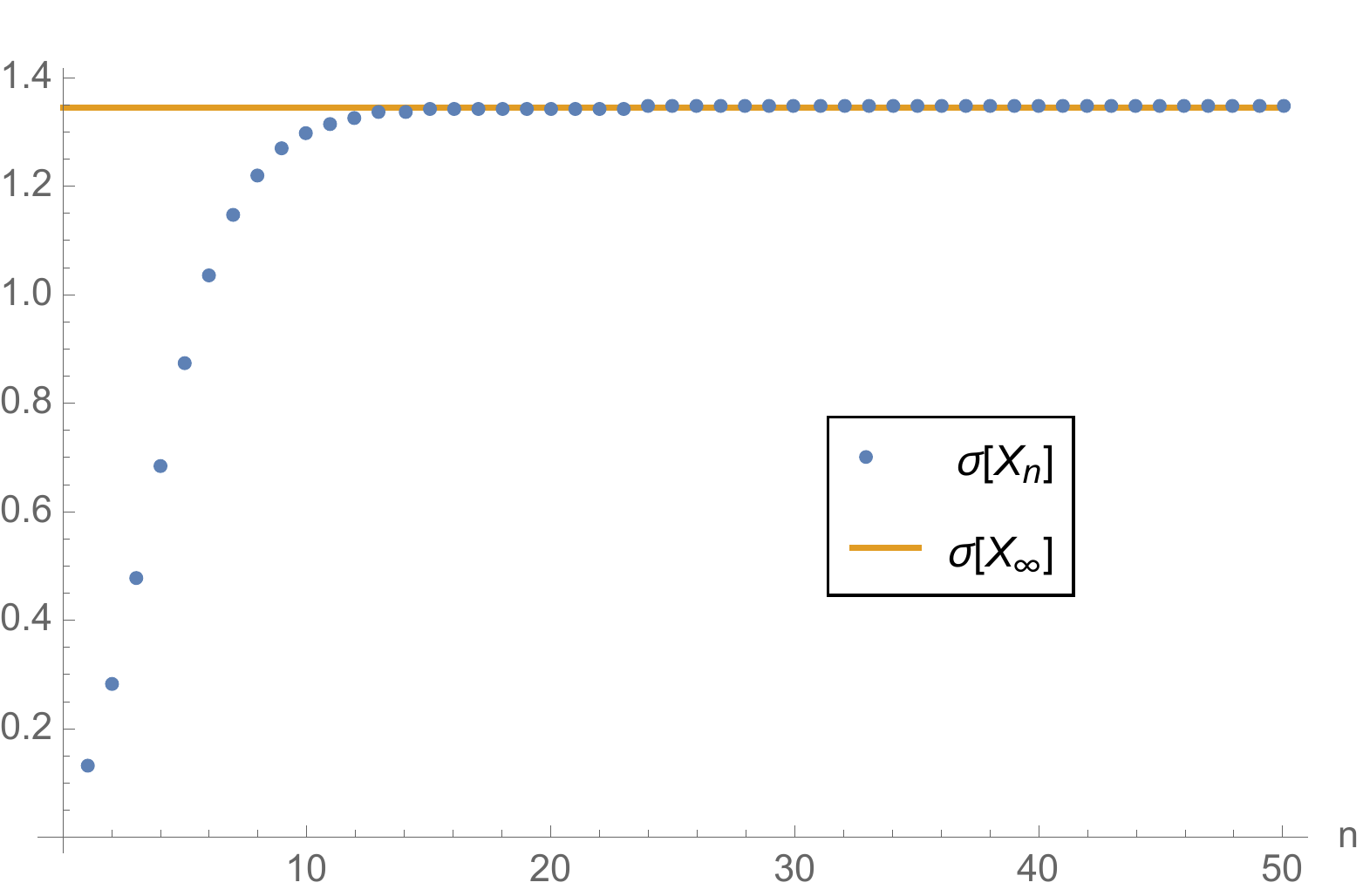}
}
\caption{Left: Blue points: $\mathbb{E}\left[X_n\right]$ for different $n\in \{1,2,\dots,50\}$. Orange solid line: $\mathbb{E}\left[X_{\infty}\right]$. Right: Blue points: $\sigma\left[X_n\right]$ for different $n\in \{1,2,\dots,50\}$. Orange solid line: $\sigma\left[X_{\infty}\right]$. Example \ref{subsec31}.}
\label{fig3}
\end{figure}

In Figure \ref{fig3} it is represented the expectation  (left) and the standard deviation (right) of $X_n$. One can observe that they tend to the expectation and the standard deviation of the steady state, respectively.

\subsection{Application: modelling the diffusion of a technology}
\label{subsec32}

This subsection is addressed to apply the theoretical results obtained in Section \ref{sec2} to model a problem using real data.
 In particular, we consider the number of mobile lines in Spain including prepaid or postpaid contracts. We assume that this situation can be modelled by the random Pielou logistic  equation  \eqref{Pielourand}. For this proposal we consider data provided by the Spanish Commission of Markets and Competition ('Comisi\'on Nacional de los Mercados y la Competencia')\cite{webtel}. In Table~\ref{data} the total number of mobile lines in Spain in each year during the range of years 1999--2015  are collected. For computations, data and outputs are expressed in ten millions of mobile lines. So, all the results presented in this section are expressed in these units.

\begin{center}
\begin{table}[htbp]
\caption{Number of mobile lines in Spain, $x_n$, at year $n$ during the range of years 1999--2015. Example \ref{subsec32}. Source: Spanish Commission of Markets and Competition ('Comisi\'on Nacional de los Mercados y la Competencia')\cite{webtel}. \label{data}}
\centering
\begin{tabular}{*{3}{c} @{\hspace{2cm}} *{3}{c} } 
{\bf Year}  & $n$ &  $x_n$   &   {\bf Year}  & $n$ &  $x_n$
\\
$1999$ & $0$  & $15\,003\,708$  &   $2008$ & $9$  & $49\,623\,339$ \\
$2000$ & $1$  & $24\,265\,059$  &   $2009$ & $10$  & $51\,052\,693$ \\
$2001$ & $2$  & $29\,655\,729$  &   $2010$ & $11$  & $51\,389\,417$ \\
$2002$ & $3$  & $33\,530\,997$  &   $2011$ & $12$  & $52\,590\,507$ \\
$2003$ & $4$  & $37\,219\,839$  &   $2012$ & $13$  & $50\,665\,099$ \\
$2004$ & $5$  & $38\,622\,582$  &   $2013$ & $14$  & $50\,158\,689$ \\
$2005$ & $6$  & $42\,693\,832$  &   $2014$ & $15$  & $50\,806\,251$ \\
$2006$ & $7$  & $45\,675\,855$  &   $2015$ & $16$  & $51\,067\,569$ \\
$2007$ & $8$  & $48\,422\,470$  &   -- & --  & -- \\
\end{tabular}
\end{table}
\end{center}

Since the model parameters $A$, $B$ and $C$ do not have a physical meaning and involve very complex factors that determine the dynamics of the number of mobile phones in Spain, we take advantage of the Central Limit Theorem in Probability to approximate them via Gaussian random variables.   
Therefore, let us take  Gaussian distributions for the random inputs $A\sim \text{N}(\mu_A;\sigma_A)$, $B\sim \text{N}(\mu_B;\sigma_B)$ and $C\sim \text{N}(\mu_C;\sigma_C)$, where the means and standard deviations are determined by adjusting the real data $x_n$ to the theoretical expectation of random Pielou logistic model. For this proposal, we find a solution of the following optimization problem that consists of minimizing the square error between real data $x_n$ and punctual predictions via the expectation ($\mathbb{E}\left[ X_n (\mu_A,\mu_B,\mu_C, \sigma_A, \sigma_B, \sigma_C) \right]$),
\[
\displaystyle 
\min_{\mu_A,\mu_B,\mu_C, \sigma_A, \sigma_B, \sigma_C}
\sum_{n=0}^{16}
\left(
x_n- \mathbb{E}\left[ X_n (\mu_A,\mu_B,\mu_C, \sigma_A, \sigma_B, \sigma_C) \right]
\right)^2
\]
where
\[
\mathbb{E}\left[ X_n (\mu_A,\mu_B,\mu_C, \sigma_A, \sigma_B, \sigma_C) \right]
=\int_{-\infty}^{\infty} 
x\, {f_1^X(x;n)}\dd\!x.
\]

Using the command \texttt{NMinimize} of  Mathematica\textsuperscript{\textregistered} software for optimization, the following parameters are obtained
\begin{equation}
\begin{array}{ccccc}
\mu_A = 1.4912, & &  \mu_B=0.095109, & & \mu_C=1.76917,\\
 \sigma_A =0.00531, & & \sigma_B=0.0025587, & &   \sigma_C= 0.0050285.
\end{array}
\label{parameters}
\end{equation}

\begin{figure}[htp]
\centerline{
		\includegraphics[width=0.5\textwidth]{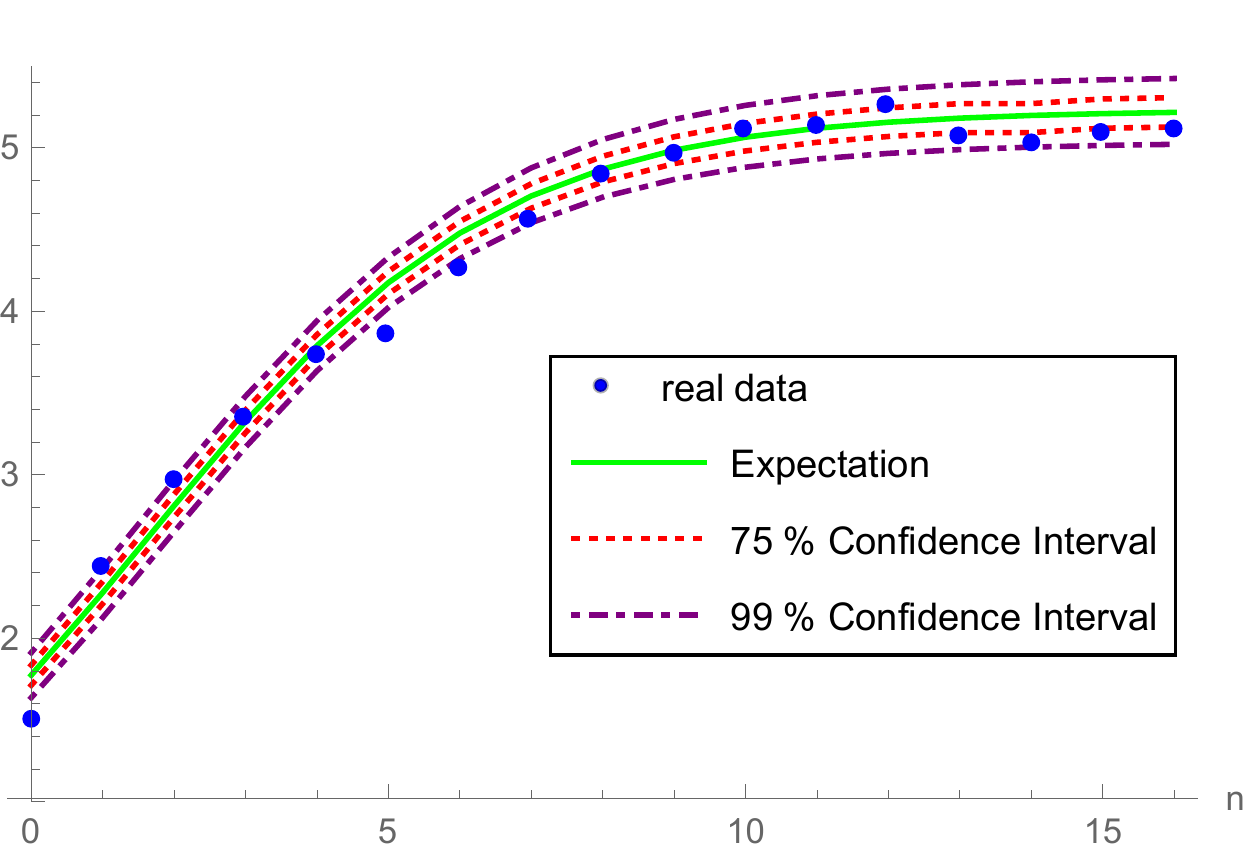}
}
	\caption{Probabilistic fitting via the proposed approach to model the dynamics of the number of mobile lines in Spain during the range of years 1999-2015 using real data collected in Table \ref{data} and $75\%$ and $99\%$ confidence intervals. Example \ref{subsec32}.}
\label{fig4}	
\end{figure}

In  Figure \ref{fig4} the number of mobile lines in Spain during the range of years 1999--2015,  $x_n$, obtained from Table \ref{data} is represented by blue points. To calculate the expectation and confidence intervals, first we calculate $f_1^X(x;n)$ with parameters \eqref{parameters}.  The expectation of $X_n$ is calculated by expression \eqref{expectation} and it has been plotted by a solid line. We can observe a good fitting between $x_n$ and the expected values $\mathbb{E}[X_n]$. In the same graphical representation, the 75\% and 99\% confidence intervals are plotted in  dash-dotted lines.  These confidence intervals have been computed in the following way. Firstly  a  value of the period $\hat{n} \geq 1$ and $\alpha\in (0,1)$ are fixed, and secondly  $z_1=z_1(\hat{n})$ and $z_2=z_2(\hat{n})$ are determined such that
\[
\int_{0}^{z_1} f_1^{X}(x;\hat{n})\,\mathrm{d}x=
\frac{z}{2}=
\int_{z_2}^{1} f_1^{X}(x;\hat{n})\,\mathrm{d}x\,.
\]
Then,  $(1-\alpha)\times 100\%$-confidence interval is specified by 
\[
1-\alpha=\mathbb{P}\left(\left\lbrace \omega\in \Omega:  
X_{\hat{n}} (\omega)\in  \left[z_1,z_2)  \right] \right\rbrace\right)
=
\int_{z_1}^{z_2} f_1^{X}(x;\hat{n})\,\mathrm{d} x\,.
\]

As we can observe most of the data lies in the confidence intervals. Only two points are outside the corresponding confidence intervals. This is a good result taking into account the uncertainties  in the market of mobile lines due to people preferences and economic oscillations.

A key feature of the previous methodology with respect to other available approaches is that it permits to construct confidence intervals considering the real probability distribution of the model output rather than using asymptotic approximations, which usually rely in the Gaussian distribution via rules of the type mean plus/minus two standard deviations.  Notice that using the  methodology proposed in this paper, and assuming appropriate distributions for the input data (in our case we have assumed that $A$, $B$ and $C$ are Gaussian), we are able to compute the theoretical probability distribution of the model output (which in our case is not Gaussian), and from it, to determine precise confidence intervals using a prefixed but arbitrary confidence level.

\begin{figure}[htp]
\centerline{
		\includegraphics[width=0.5\textwidth]{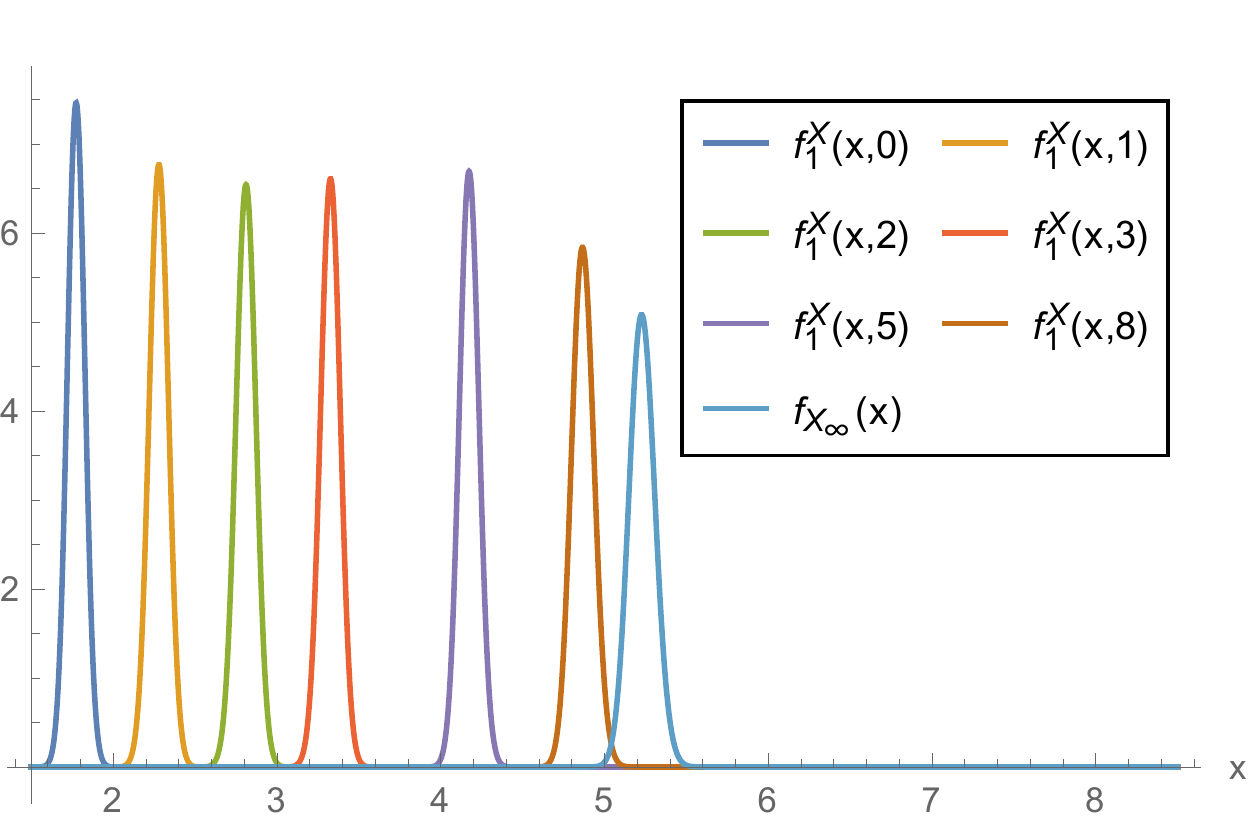}
}
	\caption{PDFs $f_1^X(x,n)$ for   different $n\in \{0,1,2,3,5,8\}$ and $f_{X_{\infty}}(x)$. Example \ref{subsec32}.}
\label{fig5}	
\end{figure}

In Figure~\ref{fig5} the 1-PDF of $X_n$, $f_1^X(x;n)$, at different fixed periods $n$ are plotted. We can observe that these 1-PDFs tend, as $n$ increases, to the PDF of the steady state, $f_{X_{\infty}}(x)$, calculated by expression \eqref{pdfinf} and   represented in light blue colour.

\begin{figure}[htp]
\centerline{
		\includegraphics[width=0.5\textwidth]{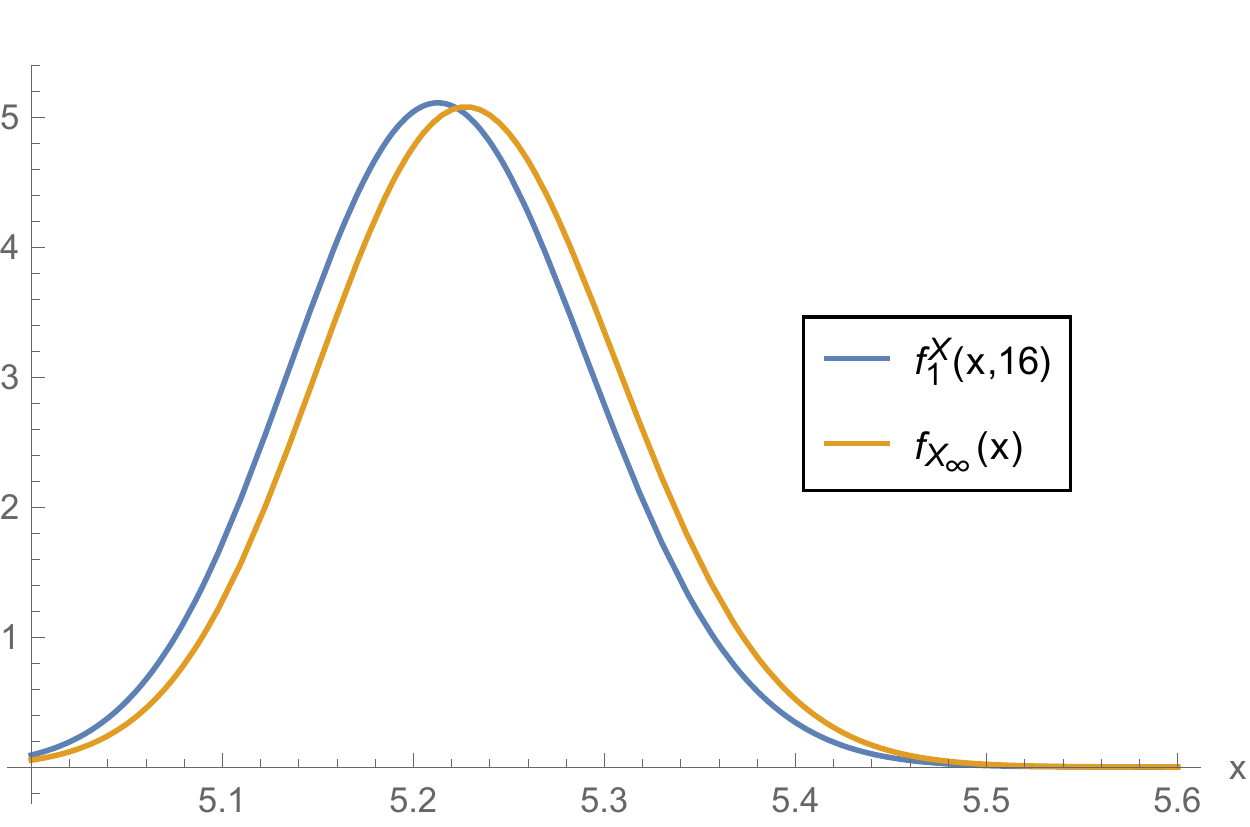}
}
	\caption{PDFs $f_1^X(x,16)$  and $f_{X_{\infty}}(x)$. Example \ref{subsec32}.}
\label{fig6}	
\end{figure}

In order to better visualize the tendency of $f_1^X(x;n)$ to $f_{X_\infty}$, as $n$ increases,  in Figure~\ref{fig6} we have plotted $f_1^X(x;16)$ together with $f_{X_{\infty}}(x)$. We can observe the mentioned tendency.

An interesting question that arises in  modelling the dynamics of the number of mobile lines is to know the period $n$ when there will be a certain number of mobile phones. This allows   to predict the period when mobile  telecommunication infrastructures might be necessary to expand, for example. The PDF of this time is calculated by \eqref{pdftiempo}.
In Figure \ref{fig6} it is plotted the distribution of time until a given proportion of the population has got  line mobile. For example, if we observe the PDF for $\hat{x}=5$, we can derive that there were about $50\,000\,000$ line mobiles in Spain at year 2008 ($n=9$).

\begin{figure}[htp]
\centerline{
		\includegraphics[width=0.5\textwidth]{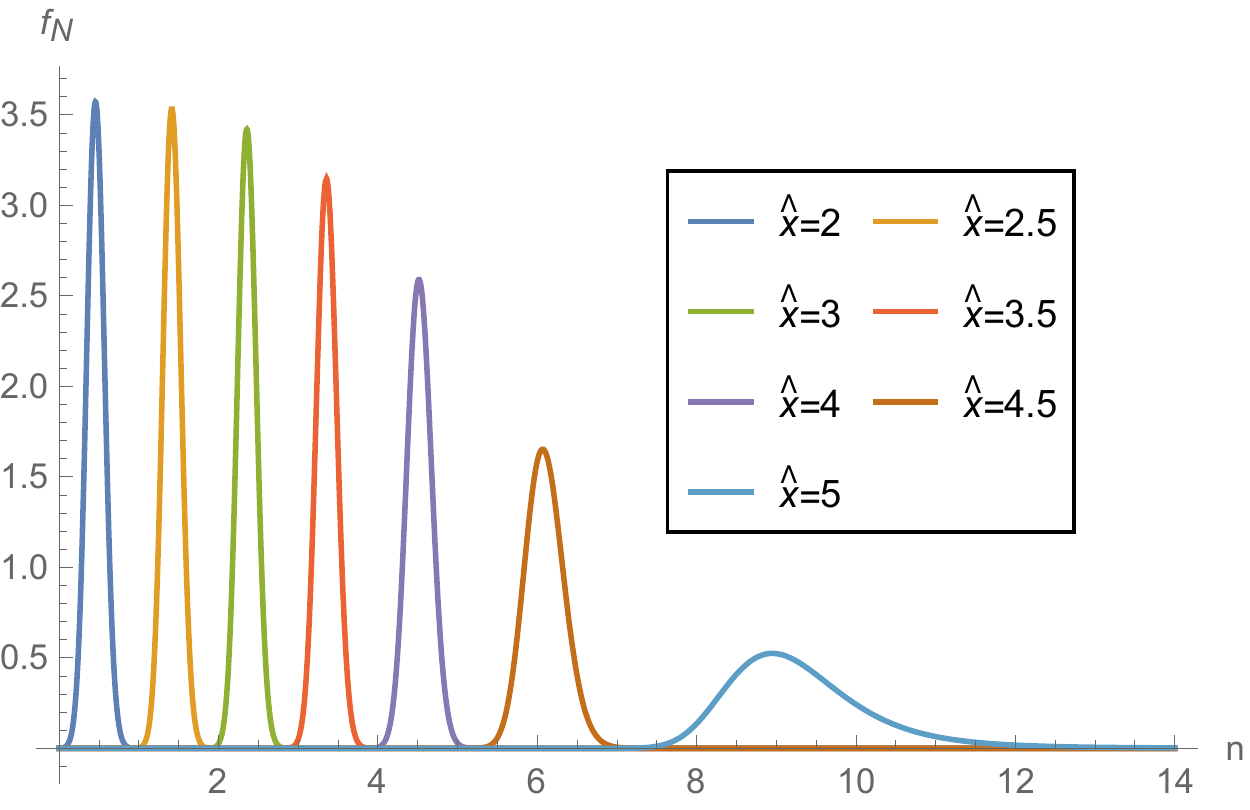}
}
	\caption{PDF of  $N$ for different sizes of the population $\hat{x}\in \{  2,2.5,3,3.5, 4,4.5,5\}$. Example \ref{subsec32}.}
\label{fig7}	
\end{figure}

\section{Conclusions}
\label{sec4}

In this paper a randomization of the Pielou logistic  model to study the dynamics of a population has been considered. 
Important goals in this type of equations are to compute its solution, to study the steady state and to compute the time until a certain level of population is reached.
A full probabilistic description of these randomized magnitudes has been provided.  This has been achieved taking advantage of the RVT technique.  

RVT technique allows us to calculate the PDF from a particular mapping.  In each case studied throughout this paper, appropriate mapping have been chosen to reduce the computational cost.

Theoretical findings have been illustrated via two  examples. The first one consists in some numerical experiments to illustrate our results from a mathematical standpoint. While the second one shows a real application to modelling the dynamics of the number of mobile lines in Spain during the range of years 1999--2015. We can conclude that the obtained results are satisfactory.

\section*{Acknowledgements}
This work has been partially supported by the Ministerio de Econom\'{i}a y Competitividad grant MTM2017-89664-P. Ana Navarro Quiles acknowledges the postdoctoral contract financed by DyCon project funding from the European Research Council (ERC) under the European Union’s Horizon 2020 research and innovation programme (grant agreement No 694126-DYCON).

\section*{Conflict of Interest Statement}
The authors declare that there is no conflict of interests regarding the publication of this article.




\section*{References}

\bibliographystyle{elsarticle-num}
\bibliography{rvtPielou}%







\end{document}